\documentclass[pre,twocolumn,showpacs,preprintnumbers,amsmath,amssymb, unsortedaddress]{revtex4}

\usepackage{dcolumn}
\usepackage{bm}
\usepackage{subfigure}
\usepackage{booktabs}
\usepackage{amsfonts,amsthm}
\usepackage{graphicx,color,amsmath,amssymb,epsfig}

\newcommand{\dt}{\mbox{d}t}
\newcommand{\du}{\mbox{d}u}

\newtheorem{theorem}{Theorem}

\begin{document}

\preprint{}

\title{Reaction Time for Trimolecular Reactions in  \\
Compartment-based Reaction-Diffusion Models}

\author{Fei Li}
\email{felix@cs.vt.edu}

\author{Minghan Chen} 
\email{cmhshirl@cs.vt.edu}

\affiliation{%
Department of Computer Science, Virginia Tech, \\
Blacksburg, VA 24061, USA
}

\author{Radek Erban}
\email{erban@maths.ox.ac.uk}
 
\affiliation{%
Mathematical Institute, University of Oxford \\
Radcliffe Observatory Quarter, Woodstock Road, 
Oxford, OX2 6GG, United Kingdom 
}%

\author{Yang Cao}
\email{ycao@cs.vt.edu}

\affiliation{%
Department of Computer Science, Virginia Tech, \\
Blacksburg, VA 24061, USA
}

\date{\today}

\begin{abstract}
\noindent
Trimolecular reaction models are investigated in the compartment-based
(lattice-based) framework for stochastic reaction-diffusion 
modelling. The formulae for the first collision time and the mean reaction
time are derived for the case where three molecules are present
in the solution. 
\end{abstract}

\maketitle

\section{Introduction}

\noindent
Trimolecular reactions are important components of biochemical models which
include oscillations~\cite{Schnakenberg:1979:SCR}, multistable 
states~\cite{Schlogl:1972:CRM} and pattern formation~\cite{Qiao:2006:SDS}. 
Considering their reactant complexes, trimolecular reactions can 
be subdivided into the following three forms 
\begin{eqnarray} 
\label{basic3} 
3 U &\displaystyle \mathop{\longrightarrow}& \emptyset, 
\\
\label{basic2} 
2 U + V &\displaystyle \mathop{\longrightarrow}& \emptyset, 
\\
\label{basic1} 
U + V + W &\displaystyle \mathop{\longrightarrow}& \emptyset, 
\end{eqnarray} 
where $U$, $V$ and $W$ are distinct chemical species (reactants)
and symbol $\emptyset$ denotes products. In what follows, we will
assume that product complexes $\emptyset$ do not include 
$U$, $V$ and $W$. Let us denote 
by $u$ the concentration of $U$. Then the conventional reaction-rate
equation for trimolecular reaction~(\ref{basic3}) can be written as
\begin{equation} 
\label{basicODE} 
\frac{\du}{\dt} = - k(t) \, u^3,
\end{equation} 
where $k(t)$ denotes (in general, time-dependent) reaction 
rate constant~\cite{Oshanin:1995:SAT}.
Using mass-action kinetics, rate $k(t)$ is assumed to be constant and equation
(\ref{basicODE}) becomes an autonomous ordinary differential equation (ODE)
with a cubic nonlinearity on its right-hand side. Cubic nonlinearities 
significantly enrich the dynamics of ODEs. For example, ODEs describing
chemical systems with two chemical species which do not include cubic
or higher nonlinearites cannot have any limit cycles~\cite{Plesa:2015:CRS}. 
On the other hand, it has been reported that, by adding cubic 
nonlinearities to such systems, one can obtain chemical systems 
undergoing homoclinic~\cite{Plesa:2015:CRS} and 
SNIC bifurcations~\cite{Erban:2009:ASC}, i.e. oscillating solutions
are present for some parameter regimes.

Motivated by the developments in systems biology~\cite{Paulsson:2000:SFF,Elowitz:2002:SGE},
there has been an increased interest in recent years in stochastic methods for
simulating chemical reaction networks. Such approaches provide 
detailed information about individual reaction events. Considering well-mixed 
reactors, this problem is well understood. The method of choice is 
the Gillespie algorithm~\cite{Gillespie:1977:ESS} or its equivalent 
formulations~\cite{Gibson:2000:EES,Cao:2004:EFS}. These methods
describe stochastic chemical reaction networks as continous-time 
discrete-space Markov chains. They are applicable to modelling 
intracellular reaction networks in relatively small domains which 
can be considered well-mixed by diffusion. If this assumption is
not satisfied, then stochastic simulation algorithms for spatially 
distributed reaction-diffusion systems have to be 
applied~\cite{Erban:2007:PGS,Erban:2007:RBC}. The most common
algorithms for spatial stochastic modelling in systems biology can 
be classified as either Brownian dynamics 
(molecular-based)~\cite{Andrews:2004:SSC,Takahashi:2010:STC} 
or compartment-based (lattice-based)
approaches~\cite{Hattne:2005:SRD,Engblom:2009:SSR}. Molecular-based 
models describe a trajectory of each molecule involved 
in a reaction network as a standard Brownian motion. This
can be justified as an approximation of interactions (non-reactive
collisions) with surrounding molecules (heat bath) on sufficiently 
long time scales~\cite{Erban:2014:MDB,Erban:2015:CAM}. It is then
often postulated that bimolecular or trimolecular reactions occur 
(with some probability) if the reactant molecules are sufficiently
close~\cite{Oshanin:1995:SAT,Erban:2009:SMR,Agbanusi:2013:CBR,Flegg:2015:SRK}.

Brownian dynamics treatment of bimolecular reactions is based on
the theory of diffusion-limited reactions, which postulates that
a bimolecular reaction occurs if two reactants are within distance
$R$ (reaction radius) from each other. The properties of this 
model depend on the physical dimension of the reactor. Considering one-dimensional~\cite{Torney:1983:DRO} 
and two-dimensional problems~\cite{Torney:1983:DRR}, diffusion-limited 
reactions lead to mean-field models with time-dependent rate constants 
which converge to zero for large times. This qualitative property is 
in one spatial dimension shared by trimolecular reactions. It has been
shown that the mean-field model~(\ref{basicODE}) can be obtained with
time-dependent rate constant $k(t)$ satisfying for large times
\begin{equation}
k(t) \approx C \sqrt{\frac{\log t}{t}},
\label{asymptoticsofk}
\end{equation}
where $C$ is a constant~\cite{benAvraham:1993:DTR,Oshanin:1995:SAT}.
Considering three-dimensional problems, the reaction radius of a 
diffusion-limited bimolecular reaction can be related with a
(time-independent) reaction rate of the corresponding mean-field
model~\cite{Erban:2009:SMR}. Trimolecular reaction (\ref{basic3}) 
can then be incorporated into three-dimensional Brownian dynamics 
simulations either directly~\cite{Flegg:2015:SRK} or as a pair of
bimolecular reactions $U+U \displaystyle \mathop{\longrightarrow} Z$
and $Z + U \displaystyle \mathop{\longrightarrow} \emptyset$, 
where $Z$ denotes a dimer of $U$~\cite{Gillespie:1992:MPI}.

Compartment-based models divide the simulation domain into 
compartments (voxels) and describe the state of the system by numbers 
of molecules in each compartment. Compartments can be both
regular (cubic lattice) or irregular (unstructured
meshes)~\cite{Engblom:2009:SSR}. Considering that the simulation
domain is divided into cubes of side $h$, the diffusive movement
of molecules of $U$ is then modelled as jumps between neighbouring
compartments with rate $d_u = D_u / h^2$, where $D_u$ is the 
diffusion constant of the chemical species $U$. In this paper, 
we will consider compartment-based stochatic reaction-diffusion 
models of the trimolecular reactions (\ref{basic3})--(\ref{basic1})
in a narrow domain $[0,L] \times [0,h] \times [0,h]$ where $L \gg h$.
Such domains are useful for modelling filopodia~\cite{Erban:2014:MSR},
but they can also be viewed as a simplification of a real three-dimensional 
geometry when there is no variation in $y$ and $z$ directions.
Then the mean-field model of trimolecular reaction~(\ref{basic3}) 
can be formulated in terms of spatially varying concentration 
$u(x,t)$, where $x \in [0,L]$ and $t \ge 0$, which satisfies 
partial differential equation (PDE)
\begin{equation} 
\label{RDPDE}
\frac{\partial u} {\partial t} = D_u \frac{\partial^2 u}{\partial x^2}  
- k \, u^3, 
\end{equation}
where $k$ is the macroscopic rate constant of the trimolecular
reaction. To formulate the compartment-based stochastic
reaction-diffusion model, we divide the domain 
$[0,L] \times [0,h] \times [0,h]$ into $K=L/h$ cubic compartments 
$[(i-1)h,ih] \times [0,h] \times [0,h] $, $i=1,2,\dots,K$. 
Denoting $U_i$ the number of molecules of $U$ in the $i$-th 
compartment, the diffusion of $U$ can be written as
the chain of ``chemical reactions"~\cite{Erban:2007:PGS}
\begin{equation}
U_1
\;
\mathop{\stackrel{\displaystyle\longrightarrow}\longleftarrow}^{d_u}_{d_u}
\;
U_2
\;
\mathop{\stackrel{\displaystyle\longrightarrow}\longleftarrow}^{d_u}_{d_u}
\;
U_3
\;
\mathop{\stackrel{\displaystyle\longrightarrow}\longleftarrow}^{d_u}_{d_u}
\;
\dots
\;
\mathop{\stackrel{\displaystyle\longrightarrow}\longleftarrow}^{d_u}_{d_u}
\;
U_K.
\label{diffGillU}
\end{equation}
Diffusion of $V$ and $W$, which appear in trimolecular reactions (\ref{basic2}) 
or (\ref{basic1}), is described using the chains 
of ``chemical reactions" of the form~(\ref{diffGillU})
with the jump rates given by $d_v = D_v / h^2$ and $d_w = D_w / h^2$,
where $D_v$ and $D_w$ are the  diffusion constant of the chemical 
species $V$ and $W$, respectively. Trimolecular chemical reactions 
are localized in compartments, i.e. each of trimolecular 
reactions~\hbox{(\ref{basic3})--(\ref{basic1})} is replaced by $K$ reactions:
\begin{eqnarray} 
\label{basic3d} 
3 U_i &\displaystyle \mathop{\longrightarrow}^{k}& \emptyset, 
\\
\label{basic2d} 
2 U_i + V_i &\displaystyle \mathop{\longrightarrow}^{k}& \emptyset, 
\\
\label{basic1d} 
U_i + V_i + W_i &\displaystyle \mathop{\longrightarrow}^{k}& \emptyset, 
\end{eqnarray} 
where $i=1,2,\dots,K$, and $k$ denotes the macroscopic reaction rate 
constant with units [m$^6$ s$^{-1}$]. Compartment-based modelling postulates 
that each compartment is well-mixed. In particular, chemical reactions
(\ref{diffGillU})--(\ref{basic1d}) can be all simulated using the 
Gillespie algorithm~\cite{Gillespie:1977:ESS} (or other algorithms 
for well-mixed chemical systems~\cite{Gibson:2000:EES,Cao:2004:EFS})
and the system can be equivalently described using the reaction-diffusion
master equation (RDME)~\cite{Erban:2007:PGS}.
In particular, the probability that a trimolecular reaction occurs
in time interval $[t,t+\Delta t)$ in a compartment containing one 
triplet of reactants is equal to $\alpha_1 \, \Delta t$ where
\begin{equation}
\alpha_1 = \frac{k}{h^6}
\label{comprate}
\end{equation}
and $\Delta t$ is chosen sufficiently small, so that 
$\alpha_1 \, \Delta t \ll 1$. The standard scaling of reaction rates 
(\ref{comprate}) is considered in this paper when we investigate
the dependence of trimolecular reactions on $h$. It has been previously
reported for bimolecular reactions that the standard RDME scaling leads
to large errors of bimolecular reactions~\cite{Erban:2009:SMR}.
One of the goals of the presented manuscript is to investigate the
effect of $h$ on trimolecular reactions.

Diffusive chain of reactions~(\ref{diffGillU}) is formulated using
a narrow three-dimensional domain $[0,L] \times [0,h] \times [0,h]$.
It can also be interpretted as a purely one-dimensional simulation.
In this case, the simulated domain is one-dimensional interval
$[0,L]$, divided into $K=L/h$ intervals $[(i-1)h,ih]$, 
$i=1,2,\dots,K$, where $U_i$ is the number of molecules of $U$ 
in the $i$-th interval. Trimolecular reactions 
(\ref{basic3d})--(\ref{basic1d}) are then written with one-dimensional
rate constant 
\begin{equation} 
	k_{1D} = \frac{k}{h^4}, 
	\label{1Drate} 
\end{equation} 
which have physical units [m$^2$ s$^{-1}$]. 
Then the probability that a trimolecular reaction occurs
in time interval $[t,t+\Delta t)$ in an interval containing one 
triplet of reactants is equal to $\alpha_2 \, \Delta t$ where
\begin{equation}
	\alpha_2 = \frac{k_{1D}}{h^2}, 
\label{comprate2}
\end{equation}
and $\Delta t$ is chosen sufficiently small, so that 
$\alpha_2 \, \Delta t \ll 1$. In particular, two interpretations of
the simulated diffusive reaction chain~(\ref{diffGillU}) considered 
in this paper differ by physical units of $k$ and $k_{1D}$, and by the corresponding
scaling of reaction rate with $h$: compare~(\ref{comprate}) and
(\ref{comprate2}). The simulated model can also be used to verify
the corresponding Brownian dynamics result~(\ref{asymptoticsofk}),
provided that we appropriately relate $h$ and $R$, and postulate 
that the trimolecular reaction occurs (for sure) whenever
a compartment contains one triplet of reactants~\cite{Oshanin:1995:SAT}.

The rest of the paper is organized as follows.
In Section \ref{secbimol}, we summarize recent results for bimolecular 
reactions in two spatial dimensions for periodic boundary conditions
and generalize them to reflective (no-flux) boundary conditions. 
These results are useful for investigating compartment-based 
stochastic reaction-diffusion modelling of trimolecular 
reactions~(\ref{basic3})--(\ref{basic1})
as we will show in Section \ref{sectrimol}. 
Computational experiments illustrating the presented results
are given in Section~\ref{numer}.
 
\section{Bimolecular reactions in two dimensions}
\label{secbimol}

\noindent
In this paper, we will focus on a special case where there is
only one combination of reactants present in the system and we
study the times for trimolecular reactions (\ref{basic3})--(\ref{basic1}) 
to fire. A similar problem for a bimolecular reaction 
$U + V \longrightarrow \emptyset$  
with one $U$ molecule and one $V$ molecule is investigated 
in Hellander et.~al.~\cite{Hellander:2012:RME} for both two-dimensional
and three-dimensional compartment-based models. 
Dividing two-dimensional domain $[0, L] \times [0, L]$ into
square compartments of size $h$, the mean time until 
the two molecules react is~\cite{Hellander:2012:RME}
\begin{equation} 
\label{taubimo} 
\tau_{bimol} =  
\frac{h^2 (1+ N^1_{\mathrm{steps}})}{k_b}
	+
\frac{h^2 N_{\mathrm{steps}}}{4 (D_u + D_v)},
\end{equation}
where $k_b$ is the bimolecular reaction rate with units [m$^2$ s$^{-1}$],
and $N^1_{\mathrm{steps}}$ is the mean number of diffusive jumps 
until $U$ and $V$ are in the same compartment, given that they 
are initially one compartment apart. $N_{\mathrm{steps}}$ 
is the mean number of steps, until $U$ diffuses 
to $V$'s location for the first time. The quantities $N_{\mathrm{steps}}$
and $N^1_{\mathrm{steps}}$ can be estimated using the theorem
in Montroll~\cite{Montroll:1969:RWL}:
\begin{theorem} \label{theorem1}
Assume that the molecule $U$ has a uniformly distributed random 
starting position on a finite two-dimensional square 
lattice with periodic boundary conditions. Then the following 
holds:
\begin{equation} \label{2Dsteps}
N_{\mathrm{steps}} = \pi^{-1} N \log(N) + 0.1951N + O(1), 
\end{equation} 
where $N = L^2/h^2$ is the number of lattice points (compartments) 
in the domain. Furthermore, $N^1_{\mathrm{steps}} = N - 1$. 
\end{theorem} 

\noindent
Using Theorem~\ref{theorem1} and (\ref{taubimo}), we
have \cite{Hellander:2012:RME} 
\begin{equation} 
\label{taubimo2}
\tau_{bimol} 
=  
\frac{L^2}{k_b}
+
\tau_{coll},
\end{equation} 
where $\tau_{coll}$ is the mean time for the first collision of
molecules $U$ and $V$, which can be approximated by
\begin{equation} 
\label{taucoll2per}
\tau_{coll}
\approx
\frac{L^2 \, \log (L \, h^{-1})}{2 \pi (D_u+D_v)} \,
+ 
4.878 \!\times\! 10^{-2}
\, \frac{L^2}{D_u+D_v},
\end{equation} 
as $h \longrightarrow 0$. Using (\ref{taubimo2})--(\ref{taucoll2per}), 
we see that the reaction time for the bimolecular reaction 
$U + V \longrightarrow \emptyset$ tends to infinity when the 
compartment size $h$ tends to zero~\cite{Hellander:2012:RME}. 
In particular, equations (\ref{taubimo2})--(\ref{taucoll2per}) imply
that the bimolecular reaction is lost from simulation
when $h$ tends to zero. Thus there has to be a lower bound 
for the compartment size. This has also been shown using 
different methods in three-dimensional systems and improvements 
of algorithms for $h$ close to the lower bound have been 
derived~\cite{Erban:2009:SMR,Hellander:2012:RME}. 

\begin{figure}[htb]
\centering
\includegraphics[width=3.0in]{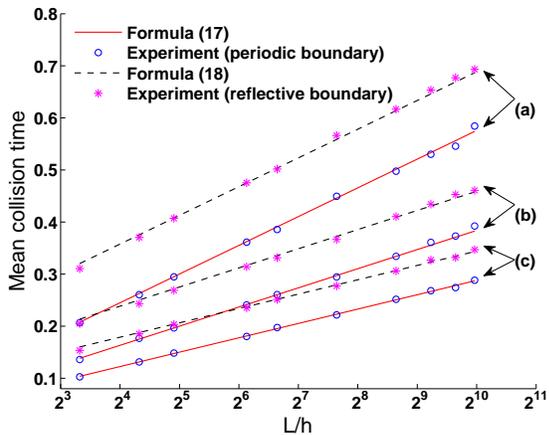}
\caption{\it Comparison of the original formula~{\rm \eqref{taucoll2per}} and our revised formula~{\rm \eqref{2DcollisionCaseA}} with mean 
collision time from numerical simulations with periodic and reflective 
boundary conditions. The two molecules diffuse freely 
in the {\rm 2D} domain, with the initial position uniformly distributed in the computational domain. Three parameter sets: (a) $D_u = 1$, $D_v = 1$; (b) $D_u = 2$, $D_v = 1$; (c) $D_u = 1$, $D_v = 3$.}
\label{fig1}
\end{figure}

\begin{figure}[htb]
	\centering
	\includegraphics[width=3.0in]{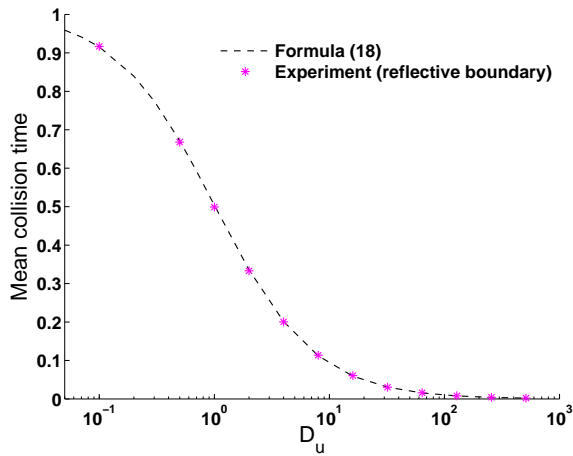}
	\caption{\it Mean collision times with reflective boundary condition 
	corresponding to different $D_u$ values when $D_v =1$ is fixed. 
	We use $L = 1.0$ and $h = 0.01$. }
	\label{fig2:caseA}
\end{figure}

Equations~(\ref{taubimo2})--(\ref{taucoll2per}) give a good approximation 
to the mean reaction time of bimolecular reaction for two-dimensional
domains, provided that the periodic boundary conditions
assumed in Theorem~\ref{theorem1} are used. However, the chain
of chemical reactions~(\ref{diffGillU}) implicitly assumes 
reflective boundary conditions. Such boundary conditions (together
with reactive boundary conditions) are commonly
used in biological systems whenever the boundary of the computational
domain corresponds to a physical boundary (e.g. cell membrane)
in the modelled system~\cite{Erban:2007:RBC}.
In Figure~\ref{fig1}, we show that formula~(\ref{taucoll2per})
is not an accurate approximation to the mean collision time 
for reflective boundary conditions.
Reflective boundary conditions mean that a molecule remains in the same 
lattice point when it hits the boundary. We plot results for
$D_u = D_v$, for $D_u > D_v$ and for $D_u < D_v$ in Figure~\ref{fig1}. 
In each case, we see that formula~(\ref{taucoll2per}) matches well with numerical results for periodic boundary conditions. In order to find a formula that matches with numerical 
experiment results for reflective boundary conditions, 
we fix the coefficient of the first term in~(\ref{taucoll2per}) and 
perform data fitting on the second coefficient. We obtain the 
following formula for the reflective boundary condition: 
\begin{equation} \label{2DcollisionCaseA} 
\tau_{coll}
\approx
\frac{L^2 \, \log (L \, h^{-1})}{2 \pi (D_u+D_v)} \,
+ 1.4053 \, \frac{L^2}{D_u+D_v}. 
\end{equation}
The average times for the first collision of $U$ and $V$
               given by (18) are plotted in Figure~\ref{fig1} corresponding to different sets of diffusion rates. 
It can be seen that formula~(\ref{2DcollisionCaseA}) matches well with numerical experiment results with reflective boundary conditions. 
In Figure~\ref{fig2:caseA}, we show another comparison between formula~(\ref{2DcollisionCaseA}) and numerical experiment data with reflective boundary conditions. 

\section {Mean Reaction Time for Trimolecular Reaction}  
\label{sectrimol} 

\noindent
In this paper, we will first consider periodic boundary conditions 
and generalize formula~(\ref{taubimo2}) 
to all cases of trimolecular reactions (\ref{basic3})--(\ref{basic1}) 
in one spatial dimension, using both scalings~(\ref{comprate}) and
(\ref{comprate2}) of reaction rates. 

Both trimolecular reactions \eqref{basic3} and \eqref{basic2} 
are special cases of~\eqref{basic1}. Since we will focus on the 
simplified situation where there is only one molecule for 
each reactant of~\eqref{basic1}, we may consider~\eqref{basic3} 
as the special case of~\eqref{basic1} where diffusion rates 
for all three molecules are the same, and~\eqref{basic2} 
as the special case where at least two of the three molecules 
have the same diffusion rates. We will denote the value of mean 
reaction time by $\tau_{trimol}$. We will decompose 
$\tau_{trimol}$ into two parts:
\begin{equation} 
\label{trimolformula}
\tau_{trimol}
=
\tau_{react}
+
\tau_{coll}, 
\end{equation}
where $\tau_{react}$ gives the mean time for reaction given
that the molecules are initially located in the same compartment,
and $\tau_{coll}$ is the mean time for the first collision, i.e.
the average time to the state of the system which has
all molecules in the same compartment,
given that they were initially uniformly distributed.
Note that the bimolecular formula (\ref{taubimo}) was
written in the form of a similar decomposition 
like~(\ref{trimolformula}). We will call $\tau_{coll}$ 
a collision time.

\subsection{Special case where $D_u = D_v $ and $D_w = 0$}
\label{firstsubDV0}

\noindent
We will start with a simple case with  $D_u = D_v $ and 
$D_w = 0$. Then the only $W$ 
molecule will be fixed at its initial location. Under 
the periodic boundary condition assumption, without loss 
of generality, we assume that this $W$ molecule is located 
at the center of the interval $[0,L]$, or specifically, 
in the compartment $\left [ \frac{K+1}{2} \right ]$, where 
$\left [ x \right ]$ represents the largest integer that 
is smaller than $x$. The $U$ and $V$ molecules diffuse 
according to (\ref{diffGillU}) from their initial 
compartments. The reaction can fire only when both 
$U$ and $V$ molecules jump to the center since the only $W$ molecule 
is at the center. Let $\tau_{coll}$ be the mean time for 
these two molecules ($U$ and $V$) to be both at the center 
for the first time, which requires that their compartment 
number is equal to $\left [ \frac{K+1}{2} \right ]$. 

Instead of trying to develop a formula for the collision time 
of two molecules diffusing to a fixed compartment in the 
one-dimensional lattice, we consider an equivalent problem: 
Imagine a $Z$ molecule jumps with a diffusion rate 
$D_u$ within an $K \times K$ grid in the two-dimensional 
space and let its compartment index be $z = (x_1, x_2)$.   
Then the two independent random walks by the two $U$ and 
$V$ molecules in one-dimensional lattice can be viewed 
equivalently as the random walk of the $Z$ molecule in 
the two-dimensional square lattice with diffusion rate $D_u$. 
Collision time $\tau_{coll}$ is then equal to the mean 
time for the pseudo molecule $Z$ to jump to the center 
for the first time, which is the case discussed in Section II. 
Therefore, the formula~\eqref{taubimo2} 
can also be applied to the trimolecular reaction~\eqref{basic1} 
when $D_w = 0$ and $D_u = D_v$ and periodic boundary conditions
are considered.  

\subsection{Special case where $D_u \neq D_v $ and $D_w = 0$}
\label{secondsubDV0}

Formula \eqref{taubimo2} cannot be directly applied to \eqref{basic1} if $D_u \neq D_v$ even 
if $D_w = 0$. We can still assume $W$ is in the center and consider the equivalent problem of a $Z$ molecule jumps with 
a diffusion rate $D_u$ in the $x$ axis and $D_v$ in the $y$ axis within an $K \times K$ grid in the 2D space. The two independent random walks by the two $U$ and $V$ molecules in 1D 
space can be viewed equivalently as the random walk of the $Z$ molecule 
in the 2D space with diffusion rate $D_u$ and $D_v$. 
We thus introduce the following theorem. The proof can be found in Appendix \ref{new2D}. 

\begin{theorem} \label{theorem2}
Assume that the molecule $Z$ has a uniformly distributed random 
starting position in a 2D lattice and that the molecules can 
move to nearest neighbours only. Assume $Z$ diffuses with diffusion rate $D_u$ 
in the $x$ direction, $D_v$ in the $y$ direction, and $D_u \geq D_v$. 
Then the following 
holds:
\begin{equation}  
\label{finalformula}
\begin{array}{rcl}
\tau_{coll} &=& \frac{\displaystyle L^2}{\displaystyle 2 \pi \sqrt{D_u D_v}} \log \left (\frac{\displaystyle L}{\displaystyle h} \right )
+  
\frac{\displaystyle L^2}{\displaystyle 12 D_v} 
\\&& +  \frac{\displaystyle L^2}{\displaystyle 4\pi  \sqrt{D_uD_v} } \left [ 2\left (\gamma + \log\left (\frac{\displaystyle 2}{\displaystyle \pi}\right) \right) - 
\log\left(1 + \frac{\displaystyle D_u}{\displaystyle D_v}\right) \right], 
\end{array}
\end{equation} 
where 
$$
\gamma = \lim \left ( 1 + \frac{1}{2} + \frac{1}{3} + \cdots 
+ \frac{1}{n} - \log n \right ) \approx 0.5772. 
$$
\end{theorem} 

\subsection{Collision time for the general trimolecular reaction}
\noindent In this subsection we consider the trimolecular reaction 
\eqref{basic1} with corresponding diffusion rates $D_u$, $D_v$ and $D_w$. 
Without loss of generality, we assume $D_u \geq D_v \geq D_w > 0$. 

\begin{figure}[htb]
\centerline{
\includegraphics[scale=0.4]{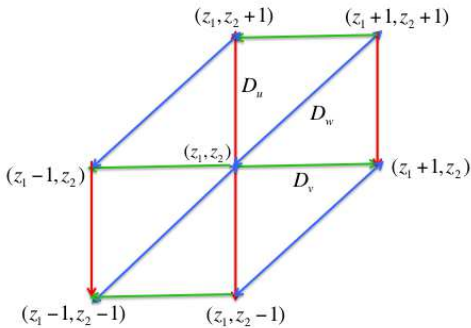}
}
\caption{{\it Illustration of the {\rm{2D}} structure for the jumps of $Z$ based on \eqref{Z2D2}. When $U$ molecule jumps, the corresponding $Z$ jumps (shown in red) the direction $(z_1, z_2) \longrightarrow (z_1, z_2 \pm 1)$. When $V$ molecule jumps, the corresponding $Z$ jumps (shown in green) the direction $(z_1, z_2) \longrightarrow (z_1 \pm 1, z_2)$; When $W$ molecule jumps, the corresponding $Z$ jumps (shown in blue) the direction 
$(z_1, z_2) \longrightarrow ( z_1 \pm 1, z_2 \pm 1)$. The whole domain has a similar shape and is not square.}}
\label{2Dillustrate2}
\end{figure}

We consider one  pseudo molecule $Z = (z_1, z_2)$, where $z_1$ and $z_2$ are expressed in terms of locations $x_u$, $x_v$
              and $x_w$ of three molecules by 
\begin{equation} \label{Z2D2}
	z_1 = {x_u - x_w}, \quad \mbox{and} 
	\quad z_2 = {x_v - x_w}. 
\end{equation} 
When $U$, $V$ and $W$ diffuses with rates $D_u$, $D_v$ and $D_w$, 
the pseudo molecule $Z$ jumps on the 2D lattice corresponding to \eqref{Z2D2}. When $Z$ jumps to the origin, $U$, $V$ and $W$ will be in the same grid, and vice versa. 
Thus the collision time 
of the trimolecular reaction in 1D is again converted to the corresponding collision time of the bimolecular problem in 2D. 
 The actual grid and jumps are illustrated in Figure \ref{2Dillustrate2}. 

\begin{figure}[htb]
\centerline{
\includegraphics[scale=0.4]{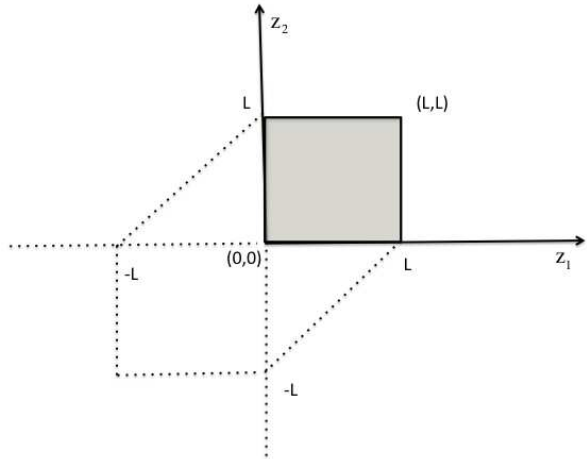}
}
\caption{{\it Illustration of the {\rm{2D}} square domain resulted from \eqref{Z2D3}. 
The original domain (dotted line) resulted from \eqref{Z2D2} is not square. With \eqref{Z2D3}, the shadowed domain is square and the theoretical analysis in Appendix can be applied.}}
\label{square}
\end{figure}

The difficulty of the mapping \eqref{Z2D2} is that the resulted domain is not square. So we cannot apply the theoretical results presented in the Appendix. 
In order to apply an estimation on a square lattice, we need to further modify the mapping \eqref{Z2D2}. 
We will take
\begin{equation} \label{Z2D3}
	z_1 = |x_u - x_w|, \quad \mbox{and} 
	\quad z_2 = |x_v - x_w|. 
\end{equation} 
Then the molecule $Z$ with coordinates $(z_1, z_2)$ will jump in an $K \times K$ 2D square lattice, where $K = L/h$. Of course, the jumps at the boundary will be different from the illustration in Figure \ref{2Dillustrate2}, 
but that is just symmetric reflection.
It does not change the validity of the formula derived in the Appendix, 
which is based on the assumption of periodic lattices. The domain resulted from the mapping \eqref{Z2D3} is shown in Figure \ref{square}. 
We have the following approximation formula for the mean jump time. 
\begin{equation}  \label{final} 
\begin{array}{rcl}
\tau_{coll} &=& \dfrac{L^2}{2 \pi \widehat D}\log\left( \dfrac{L}{h} \right)+  
\dfrac{L^2}{12 \left (D_v + D_w \right )} 
\\&& + \dfrac{L^2}{4\pi \widehat D} \left[ 2\Big(\gamma + \log\left(\dfrac{2}{\pi}\right)\Big) - 
\log\left(1 + \eta'
\right) \right]. 
\end{array}
\end{equation} 
where 
\begin{equation}  \label{hatD} 
	\widehat D = \sqrt{D_uD_v + D_uD_w + D_vD_w}, 
\end{equation} 
and
\begin{equation}  \label{etaPrime} 
	\eta' = \frac{D_u^2}{{\widehat D}^2}.  
\end{equation} 

\noindent {\bf Remark:} Formula \eqref{final} also implies an estimation 
of the collision time of bimolecular reaction in a 1D domain. 
In \eqref{final}, if we let $D_u \longrightarrow \infty$, the formula leads to 
an estimation of the mean collision time of a bimolecular reaction 
$V+W \to \emptyset$
in a 1D domain $[ 0, L]$, which is independent of $h$ 
\begin{equation} \label{Bimolecular}
	\tau_{coll} = \frac{L^2}{12 (D_v + D_w)}. 
\end{equation} 

\subsection{Formula under reflective boundary conditions}
\noindent As shown in Section \ref{numer}, formula \eqref{final} matches with 
the result computed by stochastic simulations for periodic boundary condition 
very well. However, when we consider a reflective boundary condition, 
we see a mismatch. In order to find a formula that matches with the 
reflective boundary condition results, we performed numerical experiments
and collected the mean collision time with reflective boundary conditions. 
Then with data fitting, we managed to find the formula that matches with numerical results corresponding to reflective boundary conditions. 
(See Appendix~\ref{specialcase} for  derivation of the coefficient of $0.140$.)
Equation~\eqref{collision:reflective} gives an estimation of the mean 
first collision time of 1D trimolecular reaction 
\begin{equation} \label{collision:reflective}
\begin{array}{rcl}
\tau_{coll}& = &\dfrac{L^2}{2\pi \widehat{D}} \log\left(\dfrac{L}{h}\right) + 
0.140\dfrac{L^2}{D_v+D_w} \\
&& + \dfrac{L^2}{4\pi\widehat{D}}\left[2\left(\gamma + \log\left(\dfrac{2}{\pi}\right)\right)
-\log\left(0.125+\dfrac{\eta'}{4}\right)\right],
\end{array}
\end{equation}
where $\widehat{D}$ is defined in Equation \eqref{hatD} and
$\eta^\prime$ is defined in Equation \eqref{etaPrime}.

\noindent {\bf Remark:}
When we let $D_u \rightarrow \infty$, we have 
\begin{equation} \label{Bimolecular2}
	\tau_{coll} = 0.14 \frac{L^2}{D_v + D_w}. 
\end{equation}

\subsection{Mean Reaction Time}
Formula \eqref{final} gives the estimation for the mean collision time $\tau_{coll}$ of trimolecular reactions. For the first reaction time $\tau_{trimol}$, since our derivation is based on corresponding analysis on 2D grids, an estimation similar to the equation \eqref{taubimo} can be applied. 
We have
\begin{equation} 
\label{tautrimol}
\tau_{trimol} 
\approx 
\frac{h^2 (1+ N^1_{\mathrm{steps}})}{k_{1D}}  + \tau_{coll}, 
\end{equation} 
as $h \longrightarrow 0$, where $k_{1D}$ is the reaction rate for the trimolecular reaction as defined in  \eqref{1Drate} and \eqref{comprate2}.  

For biochemical reactions with reflective boundary conditions,
we have the formula for the corresponding mean reaction time:  
\begin{equation} 
\label{tautrimol2} 
\begin{array}{rcl}
\tau_{trimol} & = &
\dfrac{L^2}{k_{1D}}  + \dfrac{L^2}{2\pi \widehat{D}} \log\left(\dfrac{L}{h}\right) + 
0.140\dfrac{L^2}{D_v+D_w} \\
&& + \dfrac{L^2}{4\pi\widehat{D}}\left[2\left(\gamma + \log\left(\dfrac{2}{\pi}\right)\right)
-\log\left(0.125+\dfrac{\eta'}{4}\right)\right]. 
\end{array}
\end{equation} 

Correspondingly, if we use the scaling in \eqref{comprate}, \eqref{tautrimol2} takes the form
\begin{equation} 
\label{tautrimol3} 
\begin{array}{rcl}
\tau_{trimol} & = &
\dfrac{L^2h^4}{k}  + \dfrac{L^2}{2\pi \widehat{D}} \log\left(\dfrac{L}{h}\right) + 
0.140\dfrac{L^2}{D_v+D_w} \\
&& + \dfrac{L^2}{4\pi\widehat{D}}\left[2\left(\gamma + \log\left(\dfrac{2}{\pi}\right)\right)
-\log\left(0.125+\dfrac{\eta'}{4}\right)\right]. 
\end{array}
\end{equation} 
We can see that as $h \rightarrow 0$, the reaction time will be dominated by the collision time in either case. 

\section{Numerical Results} \label{numer}

\begin{figure}[htb]
\centerline{
\mbox{\includegraphics[scale=0.4]{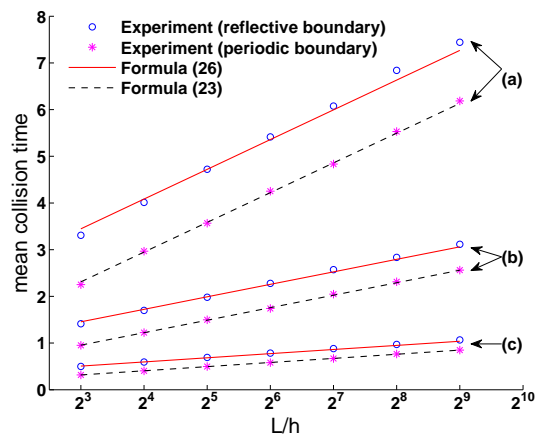}}
}
\caption{{\it The mean first collision time of three molecules in the 
{\rm{1D}} domain  (for periodic and reflective boundary conditions). 
The solid lines plot the mean colision time
with periodic and reflective boundary conditions corresponding to Formula
\eqref{final} and Formula \eqref{collision:reflective},
while dash lines show the numerical results by stochastic simulations.
For the same color from top to bottom, 
the three sets of diffusion rate parameters are
(a)
$D_u = D_v = D_w = 0.1$; 
(b) $D_u = 0.5$, $D_v = 0.2$, $D_w = 0.1$ and (c) 
$D_u = 2.5$, $D_v = 0.5$, $D_w = 0.1$.
}}
\label{collisiontime:1}
\end{figure}

\begin{figure}[htb]
\centerline{
\mbox{\includegraphics[scale=0.4]{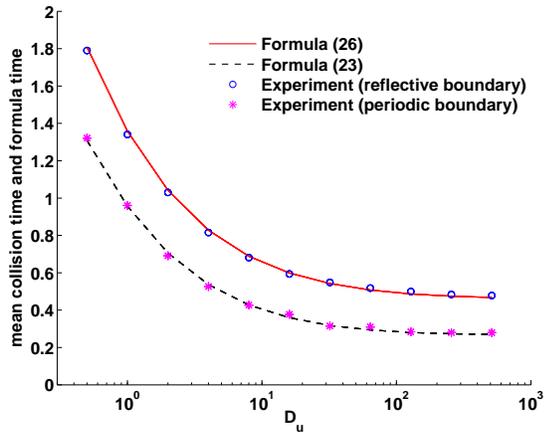}}
}
\caption{\label{collision time:2}
{\it The plot of mean trimolecular collision time (for periodic and reflective boundary conditions)
against diffusion rate $D_u$. Other parameters are given as $D_v = 0.2$, $D_w = 0.1$ and $N = 20$.
The 1D domain size is fixed at $L = 1.0$.}} 
\end{figure}

\noindent We test the collision time of the general trimolecular reaction 
\eqref{basic1}. Figure~\ref{collisiontime:1} shows the comparison of the 
numerical results of the mean first collision time 
with periodic and reflective boundary conditions with the two formulas 
\eqref{final} and \eqref{collision:reflective}. 
Figure~\ref{collisiontime:1} demonstrates that formula \eqref{final} matches 
well with computational results obtained by stochastic simulations
corresponding to periodic boundary conditions, justifying our analysis. 
It also shows that formula \eqref{collision:reflective} 
matches well with experimental results corresponding to reflective boundary conditions. 
Figure~\ref{collision time:2} shows the plot of numerical results of the 
mean first collision time with periodic and reflective boundary conditions corresponding to different $D_u$ when we fix $L$, $D_v$ and $D_w$. We can see that the two formulas match with the numerical results very well. 
We also present another comparison and analysis for the special case when $D_u \rightarrow \infty$ in Appendix \ref{specialcase}. 

\begin{figure}[htb]
\centerline{
\mbox{\includegraphics[scale=0.4]{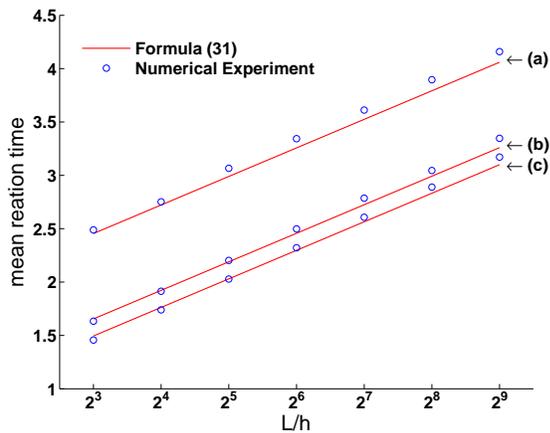}}
}
\caption{\label{reactiontime}
{\it The plot of mean trimolecular reaction time for reflective boundary 
condition against $\log(L/h)$. 
The parameters for the plot are $L=1.0, D_u=0.5, D_v=0.2, D_w=0.1$. 
For the same color, the reaction rates from top to bottom are (a) 
$k_{1D} = 1.0, (b)\ k_{1D} = 5.0, (c)\ k_{1D} = 25.0$, respectively.}}
\end{figure}

For the mean reaction time comparison, we focus only on reflective boundary condition. In Figure \ref{reactiontime}, we show the comparison of mean reaction time with formula \eqref{tautrimol2} corresponding to a fix set of diffusion rates and different reaction rates $k_{1D}$. 

\section{Discussion} 

The important consequence of Formula \eqref{final} and 
Formula \eqref{collision:reflective} 
is that when $h \longrightarrow 0$, $\tau_{coll} \longrightarrow \infty$. 
Since trimolecular reaction time $\tau_\mathrm{trimol} > \tau_{coll}$, 
the reaction time also tends to infinity.  
In particular, our analysis confirms the observations made in~\cite{Oshanin:1995:SAT, benAvraham:1993:DTR} that the diffusion-limited three-body reactions do not
recover mean-field mass-action results in 1D. The rate constant
in \eqref{basicODE}  for such problems is time dependent, given by \eqref{asymptoticsofk}, and
converges to zero as $t \to \infty$. In the limit $h \to 0$,
the average collision time goes to infinity and the rate constant
k(t) in \eqref{basicODE}  converges to zero for any finite time. It is given
by \eqref{asymptoticsofk}, where $C \equiv C(h)$ converges to zero as $h \to 0$.
In recent years, this topic has been discussed for the bimolecular reaction in the 2D and 3D cases as well. This 
limitation is a great challenge for spatial stochastic simulation. 

\medskip

\noindent{\bf Acknowledgments.}
The research leading to these results has received funding from the
European Research Council under the \textit{European Community}'s
Seventh Framework Programme ({\it FP7/2007-2013})/ERC {\it grant
agreement} n$^\circ$ 239870. Radek Erban would like to thank the Royal 
Society for a University Research Fellowship and the Leverhulme 
Trust for a Philip Leverhulme Prize. This work was carried out during the visit of Radek Erban to the
 Isaac Newton Institute. This work was partially supported by
 a grant from the Simons Foundation.

\begin{appendix}

\section{Formula of Nonuniform Random Walk on a 2D Lattice}  \label{new2D}
We consider the random walk problem on a square 2D lattice. The estimation procedure 
we present here is generalized from the original idea in Montroll~\cite{Montroll:1969:RWL}.  
Following the notations in Montroll~\cite{Montroll:1969:RWL}, let $F_n(s)$ be the probability that a lattice walker which starts at the origin arrives at a lattice point $s$ for the first time after $n$ steps and 
\begin{equation*}
	 F(s, z) = \sum_{n=1}^{\infty} z^n F_n(s) 
\end{equation*}
be the generating function of the set $\{F_n(s)\}$. Montroll showed that 
\begin{equation*} 
	F(s,z) = [P(s,z) - \delta_{s,0}]/P(0,z), 
\end{equation*} 
where $P(s, z)$ is the generating function 
\begin{equation*} 
	P(s, z) = \sum_{n=0}^{\infty} z^n P_n(s), 
\end{equation*} 
where $P_n(s)$ is the probability that a walker starting from the origin attires at $s$ for the first time after $n$ steps, no matter how many previous visits he already had at $s$. The generating function for the probability that a walker will be 
trapped at the origin in a given number of steps is 
\begin{eqnarray*} 
G_s(z) &=& \frac{1}{N-1} \sum_{s\neq 0} F(s, z) \nonumber \\
&=&  \frac{1}{N-1} \left \{ \sum_{s} F(s, z) - F(0, z) \right\}, 
\end{eqnarray*} 
where $N = m^2$ (with $m = \frac{L}{h}$) is the total number of grids 
in the square 2D lattice. 
The average number~\cite{Montroll:1969:RWL} of steps required to reach 
the origin for the first time is
\begin{eqnarray} 
	\langle n \rangle &=& \partial{G_s}/\partial{z}|_{z=1} \nonumber \\
	&=& 
	\frac{1}{N-1} 
	\frac{\partial}{\partial z} 
	\left\{ \frac{1}{(1-z)P(0, z)}\right \}_{z=1}. 
\label{formula_n}
\end{eqnarray}
In order to find the formula for $P(0.z)$, the structure function is defined as 
\begin{equation*}
 	\lambda(\theta)  = \sum_s p(s) \exp (i\theta s), 
\end{equation*} 
where $p(s)$ is the probability that at any step a random walker makes 
a displacement $s$. $\sum_s p(s) = 1$. 
In this 2D square lattice, we consider a general problem. Suppose the 
random walk can jump to left and right with 
rate $D_u/h^2$, to up and down with rate 
$D_v/h^2$, and to the diagonal directions with rate $D_w/h^2$. 
Let $\sigma_1^2 = \frac{D_u}{D_u+D_v + D_w}$, 
$\sigma_2^2 = \frac{D_v}{D_u+D_v + D_w}$, 
and $\sigma_3^2 = \frac{D_w}{D_u+D_v + D_w}$, 
We thus have 
\begin{equation*} 
	\sigma_1^2 + \sigma_2^2 + \sigma_3^2 = 1. 
\end{equation*} 
There are some special cases. 
In a uniform diffusion case, $\sigma_1 = \sigma_2 = \sigma_3= \sqrt{3}/3$. 
In the case $D_w = 0$, $\sigma_3 = 0$ and then when $\sigma_1 = \sigma_2$, that will be
the simple case discussed in Montroll~\cite{Montroll:1969:RWL}. Here we will simply derive the general formula and 
then discuss different special cases. 
For the probability $p(s)$ we have
$p(1,0) = p(-1, 0) = \sigma_1^2$/2 and $p(0, 1) = p(0, -1) = \sigma_2^2/2$, 
$p(1, 1) = p(-1, -1) = \sigma_3^2/2$
and 
\begin{equation*} 
	\lambda(\sigma)  = c_1 \sigma_1^2 + c_2 \sigma_2^2 + \sigma_3^2 (c_1c_2 + s_1s_2), 
\end{equation*} 
where $\theta = (\theta_1,\ \theta_2)$, $c_i = \cos \theta_i$ and $s_i = \sin \theta_i$.  Then 
\begin{equation*} 
	P(0, z) = m^{-2} \sum_{k_1 = 0}^{m-1}\sum_{k_2=0}^{m-1} \left[ 1 - z \lambda(2\pi k_1/m, 2\pi k_2/m) \right]^{-1}. 
\end{equation*} 
In order to obtain $\langle n \rangle$, we need to estimate $P(0,z)$. Following 
Montroll's work~\cite{Montroll:1969:RWL}, 
\begin{equation*} 
	P(0, z) = \frac{1}{m} \sum_{k_1 = 0}^{m-1} 
	\left [ 1 - c_1 \sigma_1^2 z \right]^{-1} f(z, \theta_1) , 
\end{equation*} 
where 
\begin{equation*} 
	f(z, \theta_1) = \frac{1}{m} \sum_{k_2 = 0}^{m-1} 
	\left [1 - c_2 w_1 - w_2 s_2 \right] ^{-1}, 
\end{equation*} 
with 
\begin{equation*} 
	w_1 
	= \frac{z(\sigma_2^2 + c_1\sigma_3^2}{1 - z c_1\sigma_1^2}, 
\end{equation*} 
and 
\begin{equation*} 
	w_2 
	= \frac{z(\sigma_3^2 s_1}{1 - z c_1\sigma_1^2}, 
\end{equation*} 
Following  Appendix A of Montroll~\cite{Montroll:1969:RWL}, 
we define $w_1 + iw_2 = \rho e^{i\phi}$. Thus 
\begin{equation*} 
	\rho^2 = w_1^2 + w_2^2. 
\end{equation*}
Then (equation (A7) in Montroll~\cite{Montroll:1969:RWL})
\begin{equation*} 
	f(z, \theta_1) = \frac{1}{\sqrt {1 - \rho^2}} \frac{1 - x^{2m}}{1 - 2x^m \cos m\phi + x^{2m}} , 
\end{equation*}
where $x$ is the smaller root of the equation 
\begin{equation} \label{eqx2} 
	x^2 - 2x/\rho + 1 = 0. 
\end{equation} 
Now we estimate $P(0,z)$ as 
\begin{equation*} 
P(0, z) = \frac{1}{m} \sum_{k_1=0}^{m-1} 
\left\{  \frac{ \left( 1 - z c_1\sigma_1^2\right)^{-1}}{\sqrt {1 - \rho^2}} \frac{1 - x^{2m}}{1 - 2x^m \cos m\phi + x^{2m}} 
\right \}. 
\end{equation*} 
First, when $k_1 = 0$, $c_1 = 1$, $s_1 = 0$, $w_2 = 0$ and $\phi = 0$. 
Then 
\begin{equation}
	\rho  = w_1 = z(1 - \sigma_1^2z)^{-1} (\sigma_2^2 + \sigma_3^2) =
		\frac{(1 - \sigma_1^2)z}{1 - \sigma_1^2 z}. 
\end{equation} 
The corresponding term in $P(0,z)$ is given by
\begin{equation*} 
\begin{array} {rcl}
	&&\frac{1}{m}(1-z\sigma_1^2)^{-1}  \frac{1 + x^{m}}{1 - x^m} 
	\left \{ 1 - (\frac{z(1-\sigma_1^2)}{1 - z \sigma_1^2})^2 \right \}^{-\frac{1}{2}} \\
	&= &\frac{1}{m}\frac{1 + x^{m}}{1 - x^m} 
	\left \{ (1 - z\sigma_1^2)^2 - z^2 (1 - \sigma_1^2)^2 \right \}^{-\frac{1}{2} } \\
	& = & \frac{1}{m}\frac{1 + x^{m}}{1 - x^m} 
	(1 - z)^{- \frac{1}{2}} \left \{1 + z - 2z\sigma_1^2 \right \}^{-\frac{1}{2} }
	 \\
	& = & \frac{1}{m}\frac{1 + x^{m}}{1 - x^m} 
	\frac{1}{\alpha} \left \{ 2(1 - \sigma_1^2) \right \}^{-\frac{1}{2} }
	\left [1 - \beta^2 (1 - 2\sigma_1^2)  
	\right ]^{-\frac{1}{2} }, 
	\end{array}
\end{equation*}
where $\alpha = \sqrt{1-z}$ and $\beta = \frac{\alpha}{\sqrt{2(1 - \sigma_1^2)}}$. 
Now we estimate $x$ defined as in equation \eqref{eqx2}. 
\begin{equation*} 
\begin{array}{rcl}
x & = & [1 - \sqrt{1 - \rho^2}]/\rho \\
	& = & \frac{1 - z\sigma_1^2 - \sqrt{1 - z} \sqrt{1 + z - 2\sigma_1^2z}}
	{z (1 - \sigma_1^2)}\\
	& = & 1 - 2\beta + 2\sigma_1^2 \beta^2 - 
	(2\sigma_1^2 - 3)\beta^3+ \cdots. 
\end{array}
\end{equation*} 
When $m$ is large, we have estimate 
\begin{equation*} 
	x^m = 1 - 2m\beta + 2m^2 \beta^2 - 4/3m^3 \beta^3 + \cdots. 
\end{equation*} 
Thus 
\begin{equation*} 
(1 - x^m)^{-1} =   \frac{1}{2m\beta} \left(1 + m\beta + \frac{1}{3}m^2\beta^2 + O(\beta^3) \right). 
\end{equation*} 
Then the $k_1 = 0$ term is given by 
\begin{equation*} 
	[m^2(1 - z)]^{-1} + \frac{1}{6(1 - \sigma_1^2)} + O(1-z)^{1/2}. 
\end{equation*} 
For other values of $k_1$ as $z \longrightarrow 1$, 
\begin{equation*} 
\begin{array}{rcl}
&&\left( 1 - c_1\sigma_1^2 \right )\sqrt {1 - \rho^2} \\
& = & 
\left [
(1- c_1\sigma_1^2)^2 - (\sigma_2^2  + c_1 \sigma_3^2)^2 - (\sigma_3^2 s_1)^2
\right ]^{1/2} \\
& = &  \sqrt{1- c_1} \left ( 1 - \sigma_2^4 - \sigma_3^4 - c_1 \sigma_1^4 \right )^{1/2}  \\
&= & \sqrt{1- c_1} \sqrt{2(\sigma_1^2\sigma_2^2 + \sigma_1^2\sigma_3^2 + \sigma_2^2\sigma_3^2)
+ \sigma_1^4 ( 1 - c_1)} \\
& = & 2 \sqrt{\sigma_1^2\sigma_2^2 + \sigma_1^2\sigma_3^2 + \sigma_2^2\sigma_3^2} 
\sin (\pi k_1/m) \\
&& \times \left[1 + \eta \sin^2 (\pi k_1/m) \right ]^{1/2}, 
\end{array}
\end{equation*} 
where 
\begin{equation} \label{defineEta} 
	\eta = \frac{\sigma_1^4}{\sigma_1^2\sigma_2^2 + \sigma_1^2\sigma_3^2 + \sigma_2^2\sigma_3^2}. 
\end{equation} 
Thus 
\begin{equation*} 
	P(0, z) = 
	[m^2(1 - z)]^{-1} + \frac{1}{2(1 - \sigma_1^2)} \left[ \frac{1}{3} + 
	\frac{\phi(0,1)}{r}\right] + O(1-z)^{1/2}, 
\end{equation*} 
where 
\begin{equation} \label{defineR} 
	r = \frac{\sqrt{ \sigma_1^2\sigma_2^2 + \sigma_1^2\sigma_3^2 + \sigma_2^2\sigma_3^2 }}{1 - \sigma_1^2},
\end{equation}  
and 
$\phi(0,1) = S_1 + S_2 + S_3$ with 
\begin{eqnarray*}{}
S_1 & = & 
\frac{1}{m} \sum_{k=1}^{m-1} \frac{1}{\sin(\pi k/m)},  \\
S_2 & = &  \frac{1}{m} \sum_{k=1}^{m-1} 
\frac{ [1 + \eta \sin^2 (\pi k/m)]^{-1/2} - 1 }{\sin(\pi k/m)}, \\
S_3 & = & \frac{1}{m} \sum_{k=1}^{m-1} \frac{1}{\sin(\pi k/m)} \left\{
\frac{2x^m(\cos \phi m - x^m)}{1 - 2x^m \cos \phi m +  x^{2m}} \right \} \\
& & \times 
[1 + \eta \sin^2 (\pi k/m)]^{-1/2}. 
\end{eqnarray*}
Following the Appendix A in Montroll~\cite{Montroll:1969:RWL}, 
there are estimates for $S_1$, $S_2$. 
\begin{equation*} 
	S_1 = \frac{2}{\pi}\left \{ \log m + 
[\gamma + \log (2/\pi)] - \frac{\pi^2}{72 m^2} + \cdots \right\}. 
\end{equation*} 
\begin{equation*} 
	S_2 = -\frac{1}{\pi} \log( 1 + \eta) + \frac{\eta\pi}{12m^2} + \cdots . 
\end{equation*} 
The estimation for $S_3$ is complicated. 
But fortunately 
$S_3$'s contribution to $P(0, z)$ can be considered small $r$, defined in \eqref{defineR}, is not too small. 
Note that we can obtain an estimation (see Appendix in Montroll~\cite{Montroll:1969:RWL}) 
\begin{equation*} 
	x = [1 - \sqrt{1 - \rho^2}]/\rho \approx 1 - 2rs + 2r^2s^2 + \cdots, 
\end{equation*} 
where $s = \sin \pi k/m$ and $r$ is given in \eqref{defineR}.  
Thus
\begin{equation*} 
	x^m \approx e^{-2rsk} (1 + O(1/m^2)). 
\end{equation*} 
If $r$ is close to zero, $x^m$ will be close to one. Then 
$\frac{x^m}{1-x^m}$ will be very large, and so will $S_3$. 
But if $r$ is large, $x^m$ will be close to zero (except for the term $k=0$) and $S_3$ will be small. 
In order to control the error from the $S_3$ term,  we always choose $D_u$, $D_v$ and $D_w$ such that $D_u \geq D_v \geq D_w$. In this way $r \geq \left [\frac{\sigma_1^2}{\sigma_2^2 + \sigma_3^2} 
\right ] ^{1/2} \geq \frac{\sqrt{2}}{2}$ and $S_3$ will remains relatively small, when $m$ is large. 
Thus we will simply disregard the $S_3$ term in our formula. 

To sum it up, when we choose $\sigma_1 \geq \sigma_2 \geq \sigma_3$ and considering $N = m^2$, 
we disregard $S_3$ and 
have an estimate of $P(0, z)$ as 
\begin{equation} \label{estimate} 
\begin{array}{rcl}
	P(0, z) &\approx& [N(1 - z)]^{-1} + \{ c_1 \log N + c_2 + c_3/N \} 
	\\&& + O(N^{-2}) + O((1-z)^{1/2}), 
\end{array}
\end{equation} 
where 
\begin{eqnarray*}
c_1 &=& \frac{1}{2\pi \hat \sigma},
\\ 
c_2 &=& \frac{1}{6 (1 - \sigma_1^2) } +  \frac{1}{2\pi  \hat \sigma
} 
\left [ 2(\gamma + \log(2/\pi)) - \log(1 + \eta) \right],\\ 
c_3 &=& \frac{\pi}{24  \hat \sigma} (\eta - 1/3), 
\end{eqnarray*} 
with 
$$
\hat \sigma = \sqrt{\sigma_1^2\sigma_2^2 + \sigma_1^2\sigma_3^2 + \sigma_2^2\sigma_3^2}, 
$$
$$
\gamma = \lim \left ( 1 + \frac{1}{2} + \frac{1}{3} + \cdots 
+ \frac{1}{n} - \log n \right ) \approx 0.5772, 
$$
and 
$\eta = \frac{\sigma_1^4}{\sigma_1^2\sigma_2^2 + \sigma_1^2\sigma_3^2 + \sigma_2^2\sigma_3^2}$
. 
According to \eqref{formula_n}, we have 
\begin{equation} \label{newformula_n2}
	\langle n \rangle  = c_1 N\log N + c_2 N + c_3 + O\left(\frac{1}{N}\right). 
\end{equation}

Now consider the special case $D_w = 0$. Then $\sigma_3 = 0$ and the domain is really a square lattice. In this case, we have 
$$
	\hat \sigma = \sigma_1 \sigma_2, 
$$
$$
	\eta = \frac{\sigma_1^2}{\sigma_2^2}, 
$$
and 
$$
	r = \frac{\sigma_1}{\sigma_2}. 
$$
If we select $D_u \geq D_v$ (thus $\sigma_1 \geq \sigma_2$), $r \geq 1$, the $S_3$ term will be 
relatively small. Then 
multiply \eqref{newformula_n2} with the average time for each jump 
$1/k$, where $k = 2(D_u + D_v)/h^2$, apply $N = (L/h)^2$, 
and disregard lower order terms, we obtain 
\begin{equation}  
\label{finalformula3}
\begin{array}{rcl}
\tau_{coll} &=& \frac{L^2}{2 \pi \sqrt{D_u D_v}} \log \left (\frac{L}{h} \right )
+  
\frac{L^2}{12 D_v} 
\\&& +  \frac{L^2}{4\pi  \sqrt{D_uD_v} } \left [ 2(\gamma + \log(2/\pi)) - 
\log\left(1 + \frac{D_u}{D_v}\right) \right]. 
\end{array}
\end{equation} 
If we assume further that $D_u = D_v$, the equation \eqref{finalformula3} is close to  \eqref{taubimo2} except a small difference due to the $S_3$ term. Note that in this case, the formula \eqref{finalformula3} is a rigorous estimate. 

If $D_w \neq 0$, and assume the 2D lattice is square, 
 we  let $N = \frac{L^2}{h^2}$ and multiply \eqref{newformula_n2} with $\frac{h^2}{2(D_u + D_v + D_w)}$, 
 we end up with a similar estimation 
\begin{equation}  
\label{finalformula4}
\begin{array}{rcl}
\tau_{coll} &=& \frac{L^2}{2 \pi \hat D} \log \left (\frac{L}{h} \right )
+  
\frac{L^2}{12 (D_v + D_w)} 
\\&& +  \frac{L^2}{4\pi \hat D } \left [ 2(\gamma + \log(2/\pi)) - 
\log\left(1 + \eta'
\right) \right]. 
\end{array}
\end{equation} 
where 
$$
	\hat D = \sqrt{D_uD_v + D_uD_w + D_vD_w}, 
$$
and
	$$\eta' = \frac{D_u^2}{{\hat D}^2}. $$

\section{Special case: $D_u \to \infty$} \label{specialcase} 

When the diffusion rate of $U$ approaches to infinity, the trimolecular
system becomes a bimolecular collision model of $V$ and $W$.
The formula~\eqref{collision:reflective} gives the
mean time for bi-molecular collisions, when $D_u \to \infty$: 
\begin{equation}
\lim_{D_u \to \infty}\tau_{coll} = 0.14 \frac{L^2}{D_v + D_w}.
\label{limit:Du}
\end{equation}
In this subsection, we investigate the mean bi-molecular collision time 
and derive the formula for bi-molecular collision time 
when the other two molecules have a same diffusion rate $D_v = D_w$. 

\begin{figure}[htb]
\centerline{
\mbox{\includegraphics[scale = 0.4]{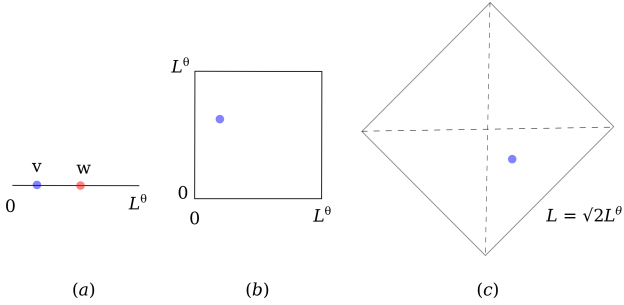}}
}
\caption{The conversion of the 1D collision model to the first exit time
model in 2D domain.
(a) The collision model of two molecules $V$ and $W$ freely diffuse in 1D.
(b) The diffusion model of one molecule in 2D domain.
(c) The first exit model of on molecule in 2D square domain.
}
\label{ftp:conversion}
\end{figure}
Assume a 1D domain of length $L^{\theta}$, two molecules $V$ and $W$ diffuse
freely with the rate $D_v$ and $D_w$ respectively.
The 1D diffusion model of two molecules is equivalent to
the 2D model in which one molecule diffuses freely with diffusion rates $D_v$ and $D_w$ 
in the two directions, independently.
The first time for the two molecules to collide in the same position 
is equivalent to the fist time when the molecule in the 2D domain comes across
 the diagonal line. 
With the reflective boundary condition, the triangle domain divided by the
diagonal line can be 
extended into a square domain with the length $L = \sqrt{2} L^\theta$. 
The first encounter time of two molecules on the 1D domain of size $L^{\theta}$
is now converted to a first exit time of on molecule on a 2D domain of size 
$L = \sqrt{2}L^{\theta}$.

In the following, we derive the formula for the first exit time
on the 2D domain. For simplicity, we assume the diffusion rate in
the two directions are the same $D = D_v = D_w$.
The Plank-Fokker equation for the diffusion in the 2D domain is
given by
\begin{equation}
\frac{\partial}{\partial t}P(\mathbf{x}, t| \mathbf{x_0}, t_0) = 
D\Big( \frac{\partial ^2}{\partial x^2} + \frac{\partial^2}{\partial y^2}\Big)
P(\mathbf{x}, t| \mathbf{x_0}, t_0), 
\label{plank_fokker}
\end{equation}
with $P(x, t|x_0, t_0)$ being the state density function, defined as the 
probability density that the molecule stays in position $\mathbf{x}$ at
time $t$ given it starts from $\mathbf{x_0}$ at time $t_0$.

Next, we define a probability function 
\begin{equation}
G(\mathbf{x_0}, t; \Omega) \equiv \int_0^L\int_0^L dxdy 
P(\mathbf{x}, t| \mathbf{x_0}, t_0),
\label{survival_prob}
\end{equation}
which describes the probability that the molecule stays in the domain 
$\Omega = (0, L) \times (0, L)$ at time $t$, 
given it starts from $\mathbf{x_0}$ at time $t=0$.
Then, we integrate the Plank-Fokker equation~\eqref{plank_fokker} over
the 2D interval $\Omega$ and we have the equation for $G$:
\begin{equation}
\frac{\partial }{\partial t} G(x_0, t; a, b) =  
D\Big( \frac{\partial ^2}{\partial x^2} + \frac{\partial^2}{\partial y^2}\Big)
G(x_0, t; a, b).
\label{pdf:G}
\end{equation}
The initial condition for the PDE~\eqref{pdf:G} is given as 
\begin{equation}
G(\mathbf{x_0}, 0; \Omega) = 1.
\label{pdf:initial}
\end{equation}
The four boundaries of the square domain $\Omega$ are all absorbing. 
Hence, we have the boundary conditions for PDE~\eqref{pdf:G} as
\begin{equation}
\begin{array}{rr}
G((0, y), t; \Omega) = 0; & \quad G((L, y), t; \Omega) = 0; \\
G((x, 0), t; \Omega) = 0; & \quad G((x, L), t; \Omega) = 0; 
\end{array}
\label{pdf:boundary}
\end{equation}

Following the definition of $G(\mathbf{x_0}, t; \Omega)$~\eqref{survival_prob},
$1-G(\mathbf{x_0}, t; \Omega)$ gives the probability that the molecule
exits $\Omega$ before time $t$, which is exactly the distribution function
of the first exit time $T(\mathbf{x_0}; \Omega)$. 
In addition, the density function of $T(\mathbf{x_0}; \Omega)$ is given by 
\begin{equation}
	p(t, \mathbf{x_0}; \Omega) 
	= \frac{\partial}{\partial t}[1-G(\mathbf{x_0}, t; \Omega)]
	= -\frac{\partial}{\partial t}G(\mathbf{x_0}, t; \Omega),
	\label{fpt:pdf}
\end{equation}
And the $n$-th moment $T_n(\mathbf{x_0}; \Omega)$ of the random variable
$T(\mathbf{x_0}; \Omega)$ is therefore given by
\begin{equation}
	T_n(\mathbf{x_0}; \Omega) 
	= \int_0^{\infty} t^n \Big[-\frac{\partial}{\partial t}
	G(\mathbf{x_0}, t; \Omega)\Big]
	dt, (n \ge 0).
	\label{fpt:moment}
\end{equation}
For $n=0$, the formula gives 
\begin{equation}
	T_0(\mathbf{x_0}; \Omega) =  - \big(G(\mathbf{x_0}, \infty; \Omega) - 
	G(\mathbf{x_0}, 0; \Omega)\big) = 1.
	\label{moment:0}
\end{equation}
Integrating~\eqref{fpt:moment} by parts, we get the formula
\begin{equation}
	\int_0^\infty t^{n-1}G(\mathbf{x_0}, t; \Omega) dt = 
	\frac{1}{n}T_n(\mathbf{x_0}; \Omega) \qquad (n\ge 1).
	\label{fpt:intpart}
\end{equation}
With equation~\eqref{fpt:pdf} and \eqref{fpt:intpart} we can formulate
the moments of $T_n(\mathbf{x_0}; \Omega)$ as coupled ordinary
differential equations. Multiply the Plank-Fokker equation~\eqref{plank_fokker}
through by $t^{n-1}$, integrating the result over all $t$ and substitute
from equation~\eqref{fpt:pdf} and \eqref{fpt:intpart}, we have the 
equation
\begin{equation}
	- T_{n-1} (\mathbf{x_0}, \Omega) = 
	D\Big(\frac{\partial^2}{\partial {x_0}^2}+\frac{\partial^2}{\partial {y_0}^2}\Big)
	(\frac{1}{n} T_n(\mathbf{x_0}; \Omega)).
	\label{moments:pde}
\end{equation}
The boundary conditions for these differential equations follows the simple
derivation from~\eqref{pdf:boundary} and we have
\begin{equation}
\begin{array}{rr}
T_n((0, y), t; \Omega) = 0; & \quad T_n((L, y), t; \Omega) = 0; \\
T_n((x, 0), t; \Omega) = 0; & \quad T_n((x, L), t; \Omega) = 0; 
\end{array}
\label{moments:boundary}
\end{equation}

With the equations ready, we can solve for the moments of the first passage
time. Here we are only interested in the first moment and
the solution to the PDE of the first moment equation yields
\begin{equation}
\begin{array}{rcl}
T_1(x, y) & = &  \displaystyle{ \frac{ x (L - x)}{2D} 
- \frac{4L^2}{D\pi^3} \sum_{k = 1, odd}^{\infty}\Big\{ 
\frac{\sin(k\pi x/L)}{k^3 \sinh(k\pi)} }\\
&&  \displaystyle{ 
\times \big(\sinh(k\pi y/L) + \sinh(k\pi (L-y)/L)\big)
\Big\} }
\end{array} 
\end{equation}
If initially the molecule is homogeneously presented in the square domain, 
we can calculate the mean first exit time as
\begin{equation}
\begin{array}{rcl}
<T_1> &=& \displaystyle{\frac{1}{L^2}\int_{0}^l \int_0^1 T_1(x, y) dxdy }\\
  & = & \displaystyle{ \frac{L^2}{12D} - \frac{16 L^2}{D}
\sum_{k=1, odd}^{\infty} \frac{\cosh(k \pi)-1}{\pi^5k^5 \sinh(k\pi)}} \\
  & \approx & 0.0351\displaystyle{\frac{L^2}{D}}
\end{array}
\end{equation}

Therefore, the mean first time when the first two molecules encounter in the 
1D domain of size $L^\theta$ is exactly the same as the mean first
exit time above and the mean first encounter time is given by 
\begin{equation}
<T_{coll}> \approx 0.0702 \frac{(L^\theta)^2}{D}.
\label{first_encounter_time}
\end{equation}

\begin{figure}[htb]
\centerline{
\mbox{\includegraphics[scale = 0.4]{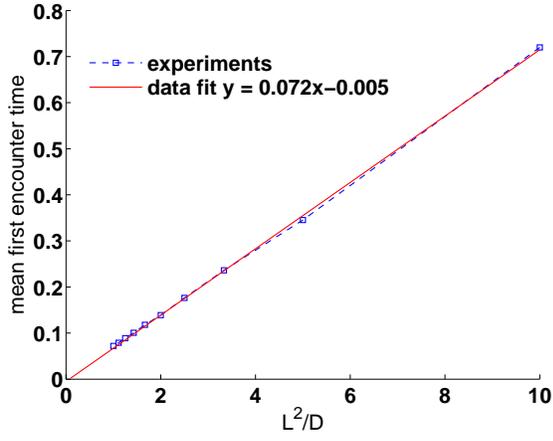}}
}
\caption{The mean first collision time of two molecules in the 1D domain with reflective boundary condition, while 
the two molecules have the same diffusion rate $D = D_v = D_w$. 
Other parameters are $N = 64$ and $L = 0.1$.
}
\label{BMEtime}
\end{figure}

\begin{figure}[htb]
\centerline{
\mbox{\includegraphics[scale = 0.4]{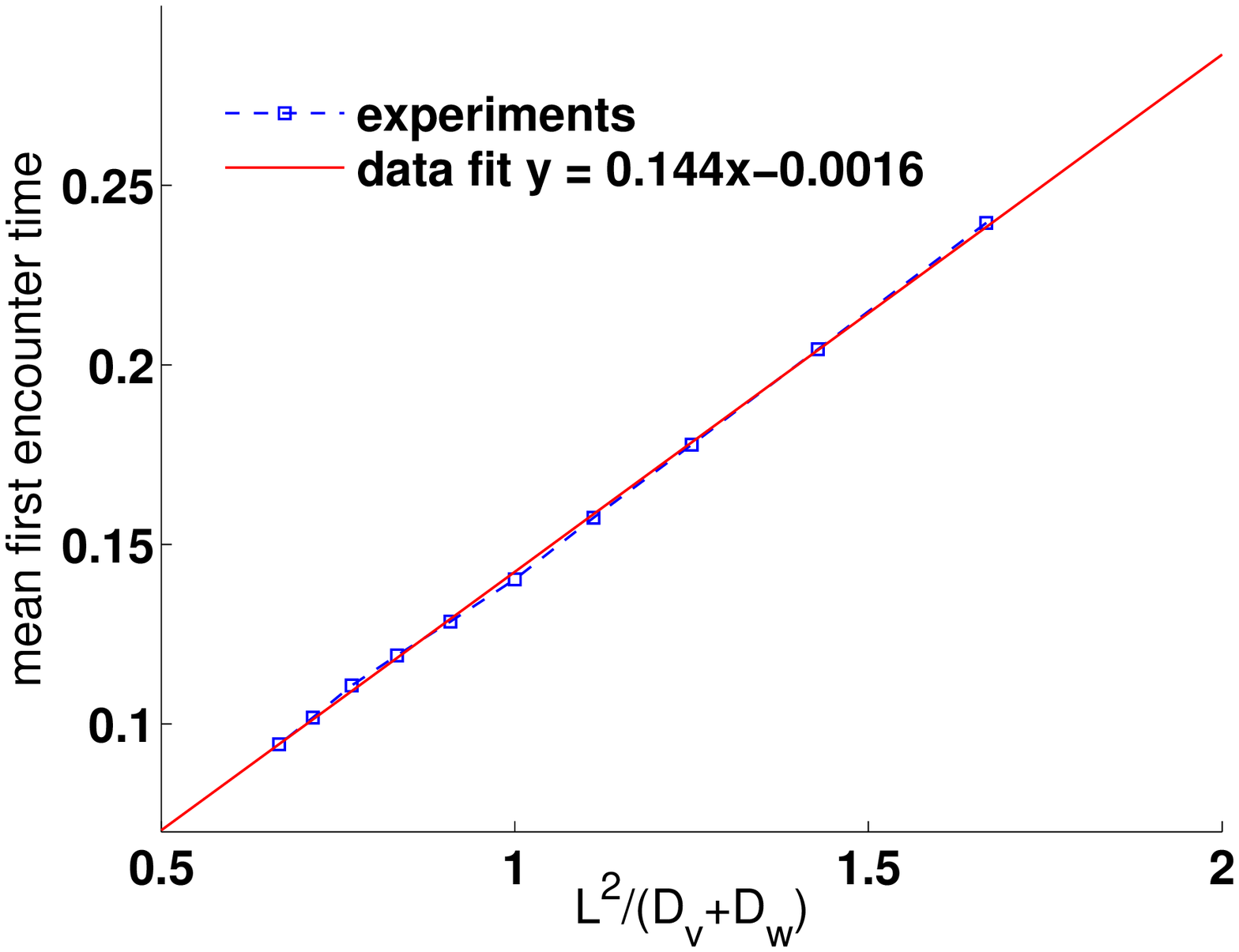}}
}
\caption{The mean first collision time of two molecules in the 1D domain with reflective boundary condition, while
the two molecules have different diffusion rates $D_v \neq D_w$. Other Parameters  are $N = 64$ and $L = 0.1$.
}
\label{BMEtime2}
\end{figure}

Figure~\ref{BMEtime} and Figure~\ref{BMEtime2} show the numerical results of the mean first 
collision time for two molecules with the same diffusion rates and for 
the general situations. The linear data fitting shows the excellent 
match with the formula~\eqref{first_encounter_time}. 
Furthermore, although our derivation is only for the case $D_v = D_w$, 
the numerical results show that the first collision time for 
the general situation, where $D_v \neq D_w$, follows the similar formula. 
This formula is given as
\begin{equation}
	<T_{coll}> \approx 0.140 \frac{(L)^2}{D_v + D_w}.
	\label{general_first_encounter_time}
\end{equation}

\begin{figure}[htb]
\centerline{
\mbox{\includegraphics[scale = 0.4]{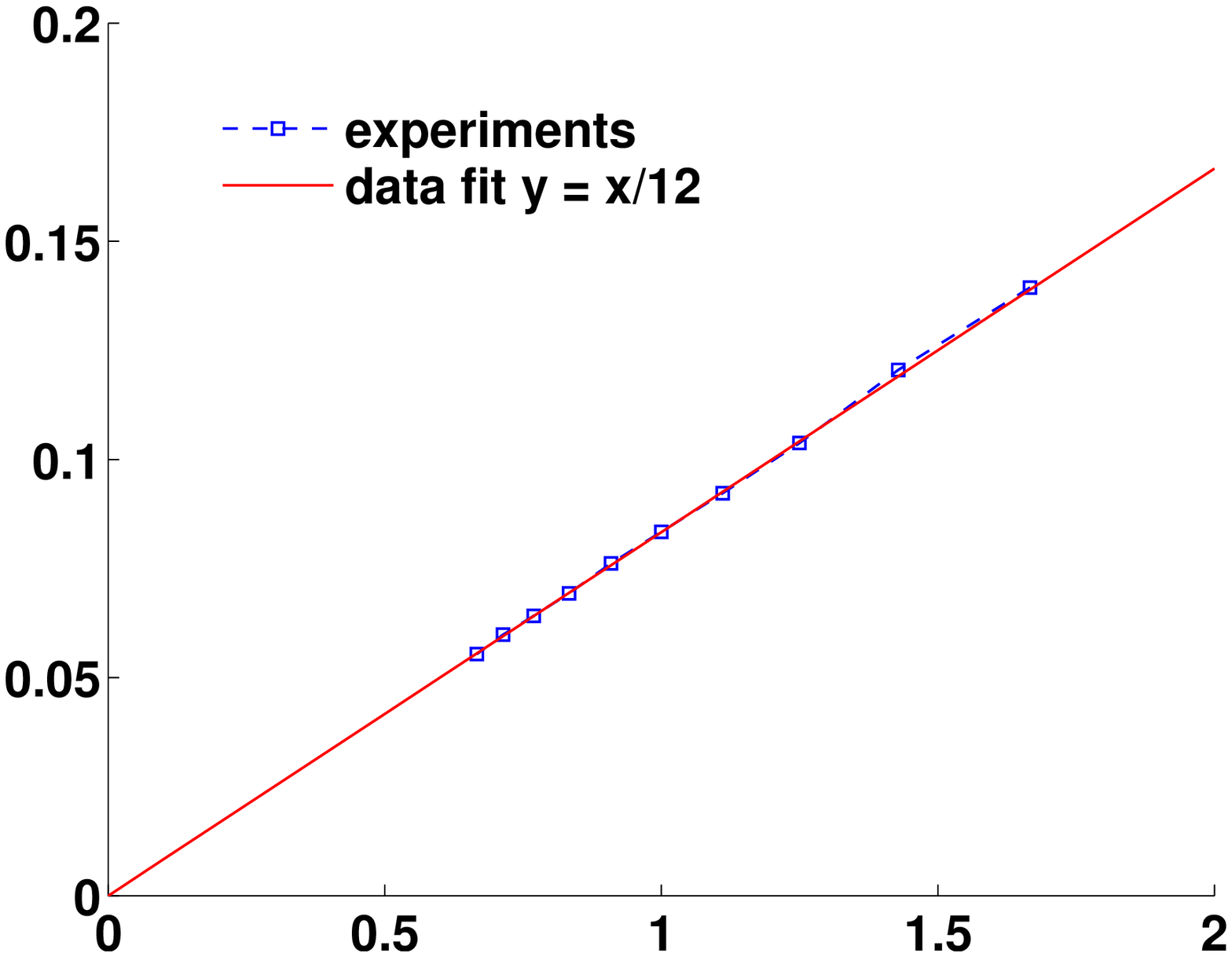}}
}
\caption{The mean first collision time of two molecules in the 1D domain with periodic boundary condition, while
the two molecules have different diffusion rates $D_v \neq D_w$. Other Parameters  are $N = 128$ and $L = 0.1$.
}
\label{periodBME}
\end{figure}

For comparison purpose, Figure \ref{periodBME} gives the numerical results of the mean first 
collision time, under periodic boundary condition,  for two molecules with different diffusion rates $D_v \neq D_w$. We can see that the numerical results match with \eqref{Bimolecular} very well. 

\end{appendix}

\end{document}